\documentclass[11pt,a4paper,reqno]{amsart}

\usepackage{amsmath}
\usepackage{amsfonts}
\usepackage{a4wide}
\usepackage{amssymb}
\usepackage{amsthm}
\usepackage{mathrsfs}
\usepackage{accents}
\usepackage[colorlinks, citecolor=blue, linkcolor=red]{hyperref}

\theoremstyle{plain}

\newtheorem{teo}{Theorem}[section]
\newtheorem{lemma}[teo]{Lemma}
\newtheorem{prop}[teo]{Proposition}
\newtheorem{cor}[teo]{Corollary}
\newtheorem{ackn}{Acknowledgments\!}

\theoremstyle{definition}

\theoremstyle{remark}

\numberwithin{equation}{section}

\def\hd{\mathring{h}}
\def\Td{\mathring{T}}
\def\SS{{{\mathbb S}}}

\def\RR{{\mathbb R}}

\def\M1{\mathscr{M}_{1}}

\def\eps{\varepsilon}

\title[Compact hypersurfaces with constant mean curvature in space forms]{A remark on compact hypersurfaces with \\constant mean curvature in space forms 
}

\author[Giovanni Catino]{Giovanni Catino}
\address[Giovanni Catino]{Dipartimento di Matematica, Politecnico di Milano, Piazza Leonardo da Vinci 32, 20133 Milano, Italy}
\email[]{giovanni.catino@polimi.it}

\begin{document}

\begin{abstract} In this note we characterize compact hypersurfaces of dimension $n\geq 2$ with constant mean curvature $H$ immersed in space forms of constant curvature and satisfying an optimal integral pinching condition: they are either totally umbilical or, when $n\geq 3$ and $H\neq 0$, they are locally contained in a rotational hypersurface. In dimension two, the integral pinching condition reduces to a topological assumption and we recover the classical Hopf-Chern result.
%
%
\end{abstract}

\maketitle

\begin{center}


\medskip

\centerline{\bf AMS subject classification:  53C40, 53C42, 53A10}

\end{center}

\

\section{Introduction}

\noindent 

The study of constant mean curvature hypersurfaces in space forms of constant curvature is one of the oldest subjects in differential geometry. There are many interesting results on this topic (see for example~\cite{hopf1, hopf2, simons, nomsmy, chedockob, chern, yau7, brendle3, andhli}, and many others). By constructing a holomorphic quadratic differential, Hopf~\cite{hopf1} showed that any constant mean curvature two-sphere in $\RR^{3}$ is totally umbilical. Chern~\cite{chern} extended Hopf's result to constant mean curvature two-spheres in three-dimensional space forms. Compact immersed constant mean curvature tori in $\RR^{3}$ were first constructed by Wente~\cite{wente}. 

To fix the notation, let $M^{n}$, $n\geq 2$, be a compact hypersurface with constant mean curvature $H$ immersed in a space form $\mathbb{F}^{n+1}(c)$ of constant curvature $c$. Denote by $h$ the second fundamental form of $M^{n}$ and by $\hd$ its trace-free part. With this notation, $M^{n}$ is totally umbilical if and only if $\hd$ vanishes. It is well known~\cite{simons, chedockob, lawson} that if $H=0$ and $|\hd|^{2} \leq nc $, $c>0$, then $M^{n}$ is either totally umbilical or a Clifford tori in $\SS^{n+1}(c)$, i.e. product of spheres $\SS^{n_{1}}(r_{1})\times\SS^{n_{2}}(r_{2})$, $n_{1}+n_{2}=n$, of appropriate radii. This rigidity result was extended by Alencar and do Carmo~\cite{aledoc} to hypersurfaces with constant mean curvature. The aim of this note is to show a characterization of compact hypersurfaces with constant mean curvature satisfying an {\em integral pinching condition} on $\hd$. This improves the result in~~\cite{aledoc}. Moreover, in dimension two, the integral inequality reduces to a topological assumption on the surface and leads to a new proof of Hopf-Chern Theorem. Our main result reads as follows:
 
\begin{teo} \label{t-main}
Let $M^{n}$ be a compact hypersurface with constant mean curvature immersed in a space form $\mathbb{F}^{n+1}(c)$ of constant curvature $c$. Then
$$
\int_{M^{n}} |\hd|^{\frac{n-2}{n}} \Big( \frac{1}{n}H^{2}-|\hd|^{2}-\frac{n-2}{\sqrt{n(n-1)}}|H||\hd|+nc \Big) \,\leq \, 0  
$$
and equality occurs if and only if $M^{n}$ is either totally umbilical or, when $n\geq 3$ and $H\neq 0$, around every non-umbilical point, it is locally contained in a rotational hypersurface of $\mathbb{F}^{n+1}(c)$.
\end{teo}

Note that, if $H=0$ and $c\leq 0$, the statement is trivial. On the other hand, there are Clifford tori in $\SS^{n+1}(c)$ with $|\hd|^{2} \equiv n c$ that are not contained in a rotational hypersurface of $\SS^{n+1}(c)$. Hence, the second part of the equality case in Theorem~\ref{t-main} cannot be true if $H=0$.

In dimension two, Gauss equation and Gauss-Bonnet theorem imply that the integral inequality is equivalent to the non-positivity of the Euler characteristic of $M^{2}$ and we recover Hopf-Chern result.

\begin{cor}[Hopf--Chern]\label{c-hopf} Let $M^{2}$ be a compact surface with constant mean curvature immersed in a space form $\mathbb{F}^{3}(c)$ of constant curvature $c$. Then, either $M^{2}$ is totally umbilical or $\chi(M^{2})\leq 0$. In particular, every  compact constant mean curvature two-sphere immersed in a space form $\mathbb{F}^{3}(c)$ is totally umbilical.
\end{cor}

The proof of Theorem~\ref{t-main} relies on an improvement of the Bochner method applied to Codazzi tensors with constant trace (section~\ref{s-cod}), which was observed by the the author in~\cite{cat3}.

\

\section{Codazzi tensors with constant trace} \label{s-cod}

Let $(M^{n},g)$ be a smooth Riemannian manifold of dimension $n\geq 3$ and consider a {\em Codazzi tensor} $T$ on $M^{n}$, i.~e., a symmetric bilinear form satisfying the Codazzi equation 
$$
(\nabla_{X} T)(Y,Z)=(\nabla_{Y} T)(X,Z)\, ,
$$ 
for every tangent vectors $X,Y,Z$. For an overview on manifolds admitting a Codazzi tensor see~\cite[Chapter 16.C]{Besse}. In all this section we will assume that $T$ has constant trace. In particular, the trace-free tensor $\Td=T-\frac{1}{n}\hbox{tr}(T)\,g$ is again a Codazzi tensor. In a local coordinate system, we have
\begin{equation}\label{eq-cod}
\nabla_{k} \Td_{ij} \,=\, \nabla_{j} \Td_{ik} \,.
\end{equation}
Throughout the article, the Einstein convention of summing over the repeated indices will be adopted. Taking the covariant derivative of the Codazzi equation and tracing we obtain
\begin{eqnarray*}
 \Delta \Td_{ij} &=& \nabla_{k} \nabla_{j} \Td_{ik} \\
 &=& \nabla_{j} \nabla_{k} \Td_{ik} - R_{ikjl} \Td_{kl} + R_{jk} \Td_{ik} \,,
\end{eqnarray*}
where we have used the commutation rules of covariant derivatives of symmetric two tensors. Here $R_{ikjl}$ and $R_{jk}$ denote the components of the Riemann and Ricci tensor respectively. Now, since $\Td$ is trace-free, from~\eqref{eq-cod} one has $\nabla_{k} \Td_{ik} = \nabla_{i} \Td_{kk} = 0$. Thus, any trace-free Codazzi tensor $\Td$ satisfies the following elliptic system
\begin{equation}\label{eq-ell}
\Delta \Td_{ij} \,=\, - R_{ikjl} \Td_{kl} + R_{jk} \Td_{ik}\,.
\end{equation}
In particular, the following Weitzenb\"ock formula holds
\begin{equation}\label{eq-wei}
\frac{1}{2}\Delta |\Td|^{2} \,=\, |\nabla T|^{2} - R_{ikjl} \Td_{ij}\Td_{kl} + R_{jk} \Td_{ij}\Td_{ik} \,.
\end{equation}

In this section we recall a vanishing theorem for Codazzi tensor with constant trace which was proved by the author in~\cite{cat3}, following the work of Gursky~\cite{gursky2} on conformal vector fields. As first observed by Bourguignon~\cite{jpb2}, trace-free Codazzi tensor satisfies the following sharp inequality.
\begin{lemma}\label{l-kat}
Let $\Td$ be a trace-free Codazzi tensor on a Riemannian manifold $(M^{n},g)$ and let $\Omega_{0}=\{p\in M^{n}:\, |\Td|(p)\neq 0\}$. Then, on $\Omega_{0}$, 
$$
|\nabla \Td|^{2} \, \geq \, \frac{n+2}{n} |\nabla |\Td| |^{2} \,.
$$
Note that, if $n=2$, then equality holds.
\end{lemma}
From the previous equation, on $\Omega_{0}$, we therefore have
\begin{equation}\label{eq-inq}
\frac{1}{2}\Delta |\Td|^{2} \,\geq\, \frac{n+2}{n}|\nabla |\Td||^{2} - R_{ikjl} \Td_{ij}\Td_{kl} + R_{jk} \Td_{ij}\Td_{ik} \,.
\end{equation}
Moreover, to apply~\eqref{eq-inq} on the whole $M^{n}$, we need to measure the set $M^{n}\setminus\Omega_{0}$. We have the following result~\cite{cat3}.

\begin{lemma}\label{l-kaz}
Let $\Td$ be a, non-trivial, trace-free Codazzi tensor on the Riemannian manifold $(M^{n},g)$ and let $\Omega_{0}=\{p\in M^{n}:\, |\Td|(p)\neq 0\}$. Then $\hbox{Vol}\,(M^{n}\setminus\Omega_{0}) = 0$. In particular~\eqref{eq-inq} holds in an $H^{1}$-sense on $M^{n}$.
\end{lemma}

Using equation~\eqref{eq-inq}, an integration by parts argument implies the following integral inequality on trace-free Codazzi tensor~\cite{cat3}.

\begin{prop}\label{p-est}
Let $\Td$ be a, non-trivial, trace-free Codazzi tensor on a compact Riemannian manifold $(M^{n},g)$. For $\eps>0$, define $\Omega_{\eps}=\{p\in M^{n}:\, |\Td|(p)\geq \eps \}$, and 
$$
f_{\eps} \,=\left\{
\begin{array}{ccc}
|\Td|(p) &\hbox{if} & p\in \Omega_{\eps} \\
\eps &\hbox{if}& p\in M^{n}\setminus\Omega_{\eps} \,.
\end{array} \right.
$$
Then
$$
\int_{M^{n}} \big(- R_{ikjl} \Td_{ij}\Td_{kl} + R_{jk} \Td_{ij}\Td_{ik}\big) \, f_{\eps}^{-\frac{n+2}{n}} \,\leq\, 0 \,.
$$
\end{prop}

\

\section{Proof of Theorem~\ref{t-main} and Corollary~\ref{c-hopf}}

Let $\mathbb{F}^{n+1}(c)$ be an $(n+1)$-dimensional smooth Riemannian manifold with constant sectional curvature $c$ and let $M^{n}$ be an $n$-dimensional compact hypersurface immersed in $\mathbb{F}^{n+1}(c)$. For any $p\in M^{n}$ we choose a local orthonormal frame $\{e_{1},\ldots,e_{n},e_{n+1}\}$ in $\mathbb{F}^{n+1}(c)$ around $p$ such that $\{e_{1},\ldots,e_{n}\}$ are tangential to $M$. Since $\mathbb{F}^{n+1}(c)$ has constant sectional curvature $c$, Codazzi and Gauss equations read (see for instance~\cite{docarmobook})
\begin{eqnarray}
\label{eq-cod}& \nabla_{k} h_{ij} - \nabla_{j} h_{ik} \,=\, 0 \,, & \\
\label{eq-gau}& R_{ikjl} \,=\, c \,(g_{ij}g_{kl}-g_{il}g_{jk}) + h_{ij}h_{kl} - h_{il}h_{jk} \,, & 
\end{eqnarray}
where $g$ denotes the induced Riemannian metric on $M^{n}$, $Rm$ its curvature tensor and $h$ the second fundamental form of $M^{n}$. In particular, tracing Gauss equation~\eqref{eq-gau}, we get
\begin{equation} \label{eq-sca}
R \,=\, n(n-1)c + H^{2} - |h|^{2} \,,
\end{equation} 
were $R$ and $H$ denote the scalar curvature of $g$ and the mean curvature of $M^{n}$, respectively.

Now, if $M^{n}$ has constant mean curvature $H$, then by Codazzi equation~\eqref{eq-cod} the tensor $\hd = h - \frac{1}{n}H  g$ is a trace-free Codazzi tensor. Thus, if $\hd$ is not identically zero, namely if $M^{n}$ is not totally umbilical, then Proposition~\ref{p-est} applies and we obtain the following integral inequality
\begin{equation}\label{eq-est1}
\int_{M^{n}} \big(- R_{ikjl} \hd_{ij}\hd_{kl} + R_{jk} \hd_{ij}\hd_{ik}\big) \, f_{\eps}^{-\frac{n+2}{n}} \,\leq\, 0 \,,
\end{equation}
where 
$$
f_{\eps} \,=\left\{
\begin{array}{ccc}
|\hd|(p) &\hbox{if} & p\in \Omega_{\eps} \\
\eps &\hbox{if}& p\in M^{n}\setminus\Omega_{\eps} \,.
\end{array} \right.
$$
and $\Omega_{\eps}=\{p\in M^{n}:\, |\hd|(p)\geq \eps \}$. Using Gauss equation~\eqref{eq-gau}, a simple calculation shows
$$
- R_{ikjl} \hd_{ij}\hd_{kl} + R_{jk} \hd_{ij}\hd_{ik} \,=\, \frac{1}{n}H^{2}|\hd|^{2}-|\hd|^{4} - H\, \hd_{ij}\hd_{ik}\hd_{jk}+nc|\hd|^{2} \,. 
$$
Moreover, since $\hd$ is trace-free, we have the sharp Okumura inequality (for a proof, see for instance~\cite[Lemma 2.6]{aledoc})
\begin{equation}\label{eq-oku}
\hd_{ij}\hd_{ik}\hd_{jk} \,\geq\, - \frac{n-2}{\sqrt{n(n-1)}}|\hd|^{3}  
\end{equation}
and, if $n\geq 3$, equality occurs at some point $p\in M^{n}$ if and only if $\hd$  can be diagonalized at $p$ with $(n-1)$-eigenvalues equal to $\lambda$ and one eigenvalue equals to $-(n-1)\lambda$, for some $\lambda\in\RR$. Hence, we obtain
$$
- R_{ikjl} \hd_{ij}\hd_{kl} + R_{jk} \hd_{ij}\hd_{ik} \,\geq\, |\hd|^{2} \Big( \frac{1}{n} H^{2}-|\hd|^{2} - \frac{n-2}{\sqrt{n(n-1)}}|H||\hd| +nc \Big)\,, 
$$
and from~\eqref{eq-est1}, we get
\begin{eqnarray*}
\int_{M^{n}}  |\hd|^{\frac{n-2}{n}} \Big( \frac{1}{n} H^{2}-|\hd|^{2} - \frac{n-2}{\sqrt{n(n-1)}}|H||\hd| +nc \Big)|\hd|^{\frac{n+2}{n}}f_{\eps}^{-\frac{n+2}{n}} \\
\leq  \int_{M^{n}} \big(- R_{ikjl} \hd_{ij}\hd_{kl} + R_{jk} \hd_{ij}\hd_{ik}\big) \, f_{\eps}^{-\frac{n+2}{n}} \leq  0 \,.
\end{eqnarray*}
Taking the limit as $\eps \rightarrow 0$, since $|\hd|^{\frac{n+2}{n}}f_{\eps}^{-\frac{n+2}{n}} \rightarrow 1$ a.e. on $M^{n}$ by Lemma~\ref{l-kaz}, we conclude
\begin{equation}\label{eq-est2}
\int_{M^{n}}  |\hd|^{\frac{n-2}{n}} \Big( \frac{1}{n} H^{2}-|\hd|^{2} - \frac{n-2}{\sqrt{n(n-1)}}|H||\hd| +nc \Big) \,\leq \, 0 
\end{equation}
and, if $n\geq 3$ and $H\neq 0$, equality occurs if and only if, at every point, either $\hd$ is null or it has an eigenvalue of multiplicity $(n-1)$ and another of multiplicity $1$. 

Now we can conclude the proof of Theorem~\ref{t-main} and Corollary~\ref{c-hopf}. If $n=2$, then we have showed that either $M^{2}$ is totally umbilic, or the following integral pinching inequality holds
$$
\int_{M^{2}}\Big(\frac{1}{2}H^{2} - |\hd|^{2} + 2c \Big) \,\leq\, 0 \,.
$$  
Since $|\hd|^{2} = |h|^{2} - \frac{1}{2} H^{2}$, we have
$$
\int_{M^{2}}\big( H^{2}-|h|^{2}+2c \big) \,\leq\, 0
$$
and from Gauss equation~\eqref{eq-sca}, we obtain
$$
\int_{M^{2}} R \,\leq \, 0.
$$
Corollary~\ref{c-hopf} now simply follows from Gauss-Bonnet theorem.

If $n\geq 3$, we have that inequality~\eqref{eq-est2} holds and equality occurs if and only if either $\hd$ is null and $M^{n}$ is totally umbilical or, when $H \neq 0$, around every non-umbilical point, $\hd$ splits with an eigenvalue of multiplicity $(n-1)$ and another of multiplicity $1$. Notice that, from Lemma~\ref{l-kaz}, the open set of non-umbilical points is dense in $M^{n}$. Theorem~\ref{t-main} now follows from~\cite[Theorem 4.2]{docdaj}, where the authors showed that every hypersurfaces in a space form with this property is contained in a rotational hypersurface of $\mathbb{F}^{n+1}(c)$.

\

\

\begin{ackn} The author is members of the Gruppo Nazionale per
l'Analisi Matematica, la Probabilit\`{a} e le loro Applicazioni (GNAMPA) of the Istituto Nazionale di Alta Matematica (INdAM) and is supported by the GNAMPA project ``Equazioni di evoluzione geometriche e strutture di tipo Einstein''. 
\end{ackn}

\

\bibliographystyle{amsplain}
\bibliography{biblio}

\

\parindent=0pt

\

\end{document}